\documentclass{commat}

\title{%
    On generalized derivations of polynomial vector fields Lie algebras
}

\author{%
    Princy Randriambololondrantomalala, Sania Asif
}

\affiliation{
    \address{D\'epartement de Math\'ematiques et Informatique, Universit\'e d'Antananarivo, Madagascar.
}
    \email{%
    princypcpc@yahoo.fr
}
    \address{School of Mathematics and Statistics, Nanjing University of information science and technology, Nanjing, Jianngsu Province, PR China.
}
    \email{%
    11835037@zju.edu.cn, 200036@nuist.edu.cn
}
}
		\thanks{We would think Pr. Ivan Kaygorodov of Federal University of ABC, Sao Paolo, Brazil, in his interest about this study. Many thanks to the referee for the valuable suggestions in this paper.}

\abstract{%
    In this paper, we study the generalized derivation of a Lie subalgebra of the Lie algebra of polynomial vector fields on $\mathbb{R}^n$, containing all constant vector fields and the Euler vector field, under some conditions on this Lie subalgebra.
}



\VOLUME{30}
\NUMBER{2}
\firstpage{51}
\DOI{https://doi.org/10.46298/cm.10386}

\begin{paper}

\section {Introduction}In~\cite{sam}, we have studied the derivations of a Lie subalgebra $ \mathfrak {P} $ of the Lie algebra of vector fields $\chi\left(\mathbb{R}^n \right)$ on $\mathbb{R}^n$ for $n\geq1$. The subalgebra $\mathfrak{P}$ contains all constant vector fields and the Euler vector field. In~\cite{pri}, we explore that each $m$-derivation $D$ of $ \mathfrak {P} $ with $m\geq2$ is an endomorphism of $ \mathfrak {P} $ such that $\forall X_1, X_2, \dots, X_m\in\mathfrak{P},$ we have
\begin{align}\label{zatra}
D&\left[X_1,\left[X_2,\dots,\left[X_{m-1},X_{m}\right]\dots\right]\right]\\
&= \left[D\left(X_1\right),\left[X_2,\dots\left[X_{m-1},X_{m}\right]\dots\right]\right]+
\left[X_1,\left[D\left(X_2\right),\dots,\left[X_{m-1},X_{m}\right]\dots\right]\right]\\
&{\quad}+\dots+\left[X_1,\left[X_2,\dots,\left[D\left(X_{m-1}\right),X_{m}\right]\dots\right]\right] +\left[X_1,\left[X_2,\dots,\left[X_{m-1},D\left(X_{m}\right)\right]\dots\right]\right].
\end{align}

It's clear that every derivation (i.e 2-derivation) of $ \mathfrak {P} $ is a $m$-derivation where $m\geq2$ but the converse is not true in general. We say that $m$-derivations for $(m\geq2)$ are a generalization of derivations. Now, we study another generalization of derivations of $ \mathfrak {P} $, that are called as generalized derivations in~\cite{luk}. However, in~\cite{luk}, only study for the finite dimensional Lie algebra has made, but study to the generalized derivations of finite and infinite dimensional Lie algebras can be made on Lie algebra $\mathfrak{P}$. We prove this fact in the form of examples for finite and infinite dimensional $\mathfrak{P}$, in the later sections. Let us recall that a \emph{monomial vector field} is a vector field on $\mathbb{R}^n$ having one and only one monomial component in one and only one element of the canonical basis of $\chi\left(\mathbb{R}^n \right)$. A Lie algebra of polynomial vector fields $ \mathfrak {P} $ is called \emph{separated} if $ \mathfrak {P} $ is spanned by monomial vector fields. Moreover, a \emph{diagonal linear vector fields} is an element of the Cartan subalgebra $\left\langle x^1 \frac{\partial}{\partial x^1},\dots,x^n \frac{\partial}{\partial x^n}\right\rangle_\mathbb{R}$ where $\left(x^i \right)_{i = 1,\dots,n}$ is the usual coordinate system of $\mathbb{R}^n$. In this research we show that if $ \mathfrak {P} $ is separated and contains all diagonal linear vector fields with some additional condition, then there are some examples of further Lie algebras structures that admit generalized derivations that were not shown in~\cite{luk} even for the finite dimensional Lie algebras. The results present in this paper are explicit and advance in the sense that we completely describe the centroid, the quasicentroid and the quasiderivations that arise in the study of generalized derivations of Lie algebras. We adapt our results in~\cite{pri} and~\cite{sam}, to study generalized derivations of infinite dimensional Lie algebras. We do not only present abstract theorems but we also give examples to illustrate all theorems of this paper.\par This paper consists of four sections after the present section for introduction. The second section is a general result from~\cite{sam} on $ \mathfrak {P} $ about its graded algebra structure and its first space of cohomology of Chevalley-Eilenberg. The third section focuses on the study of centroid and quasicentroid of $ \mathfrak {P} $ and proves that these two sets are equal. The next section describes the quasiderivation of $ \mathfrak {P} $ using the fact that this quasiderivation is a sum of homogeneous degree $i\in\mathbb{Z}$ quasiderivations of $ \mathfrak {P} $. The final section uses the results from the previous sections and presents the main theorem of this paper. According to which, generalized derivations of $ \mathfrak {P} $ are sum of derivations, quasicentroid elements, and two special endomorphisms of $ \mathfrak {P} $.
As stated earlier, it would be of great interest to explore generalized derivations of more general subalgebra $ \mathfrak {P}.$

Throughout this paper all Lie algebras are equipped as algebras by the usual bracket notation $[,]$.
\section {The Lie algebra $ \mathfrak {P} $}
Let us denote the Lie algebra of vector fields on $\mathbb{R}^{n}$ with the coordinates system $\left(x^i \right)_i$, by $\chi\left(\mathbb{R}^{n}\right)$. Throughout this section, we consider a Lie algebra $ \mathfrak {P} $ of polynomial vector fields on $\mathbb{R}^{n}.$ This Lie algebra $\mathfrak {P} $ contains the \emph{constant fields} $\frac{\partial}{\partial x^i}$ for all $i$ and the \emph{Euler vector field} $ E = {\sum_i} x ^ i\frac{\partial} {\partial x ^{i}}.$ Let $ \mathfrak{N} = \{X\in\chi \left(\mathbb{R} ^ n\right)/\left[X,\mathfrak {P}\right]\subset \mathfrak {P}\}$ denote the \emph{normalizer} of $ \mathfrak {P} $ in $ \chi \left(\mathbb{R} ^ n\right)$ and let $ \mathfrak{H}_i$ be the vector space of \emph{homogeneous polynomial vector fields} of degree $ i \in \mathbb {N}\cup\{-1\} $. Thus $X\in \mathfrak{H}_i$ if and only if
\begin{equation}
\left[E, X\right] = iX.\label{per}
\end{equation}
Denote $ \mathfrak{H} ^ {d} _ {0} $ the Lie algebra generated by all diagonal linear fields. A \emph{diagonal linear field} is identified by the diagonal matrix with coefficients in $\mathbb{R}$ of order $n$. We adopt the Einstein convention on the summation index unless expressly mentioned. Here are two examples of $ \mathfrak {P} $ in $\chi\left(\mathbb{R}^2 \right)$ with coordinate system $(x,y)$:
\begin{enumerate}
	\item The Lie algebra spanned by $\frac{\partial}{\partial x},\frac{\partial}{\partial y},x\frac{\partial}{\partial x}+y\frac{\partial}{\partial y},x^2 \frac{\partial}{\partial x}$ is of \emph{finite dimension} and
	\item The Lie algebra spanned by $\frac{\partial}{\partial x},\frac{\partial}{\partial y},x\frac{\partial}{\partial x}+y\frac{\partial}{\partial y},x^2 \frac{\partial}{\partial x},x^3 \frac{\partial}{\partial x}$ is of \emph{infinite dimension.}
\end{enumerate}
The bracket of two vector fields $X = X^i \frac{\partial}{\partial x^{i}}$ and $Y = Y^j \frac{\partial}{\partial x^{j}}$ of $\chi\left(\mathbb{R}^n \right)$ is the Lie algebra of vector fields on $\mathbb{R}^n,$ given by:
\begin{equation}
\left[X,Y\right] = \left(X^i \frac{\partial Y^j}{\partial x^{i}}-Y^i \frac{\partial X^j}{\partial x^{i}}\right)\frac{\partial}{\partial x^{j}} \label{vololona}
\end{equation}
The following lemma is clear:
\begin{lemma}\label{chris}A Lie algebra of polynomial vector fields
	$ \mathfrak {P} = \bigoplus_{i\geq -1}\mathfrak{P}_i $ is graded with the following $\left[\mathfrak{P}_i,\mathfrak{P}_j\right]\subset\mathfrak{P}_{i+j},$ where $\mathfrak{P}_i = \mathfrak{P}\cap \mathfrak{H}_i$.
\end{lemma}
We recall two theorems of~\cite{sam} on the centralizer and the cohomology space of $\mathfrak{P}$:
\begin{proposition}
Let $ \mathfrak{C} \left(\mathfrak{P}\right)$ be the centralizer of $ \mathfrak{P} $, $ \mathfrak{C} \left(\mathfrak{P}\right) = \left\{X\in \chi\left(\mathbb{R}^{n}\right)/\left[X,\mathfrak{P}\right] = \{0\}\right\}$. Here $ \mathfrak{C} \left(\mathfrak{P}\right) $ is zero. \label{lova}
\end{proposition}
\begin{definition}
We denote $Der\left(\mathfrak{P}\right)$ the Lie algebra of derivations of $\mathfrak{P}$, $L_X = \left[X,.\right]$ is the Lie derivative with respect to $X\in\mathfrak{P}$, it is called inner derivation. $ad_\mathfrak{P}$ is the Lie algebra of inner derivations of $\mathfrak{P}$. The first Chevalley-Eilenberg cohomology space is denoted $\mathrm {H} ^{1} \left(\mathfrak{P}\right) = Der\left(\mathfrak{P}\right)/ad_\mathfrak{P}$.
\end{definition}
\begin{theorem}
	$ \mathrm {H} ^{1} \left(\mathfrak{P}\right) \cong \mathfrak{N} / \mathfrak{P} \cong \mathfrak{M} $ with $ \mathfrak {M} $ a Lie subalgebra of $ \mathfrak{H}_0 $. If all linear diagonal fields are in $ \mathfrak {P} $ then $ \mathrm{H} ^{1} \left(\mathfrak{P}\right) $ is zero \label{julie}.
\end{theorem}
\section{Centroid and quasicentroid of $\mathfrak{P}$}
\begin{definition}
The centroid $C\left(L\right)$ of a Lie algebra $L$ is the Lie algebra of endomorphisms $f$ in $L$ such that
\begin{equation}
f\left[X,Y\right] = \left[fX,Y\right] = \left[X,fY\right],\label{qs}
\end{equation}
\end{definition}
\begin{definition}
The quasicentroid $QC\left(L\right)$ of a Lie algebra $L$ is the Lie algebra of endomorphism $g$ in $L$ such that
\begin{equation}
\left[g\left(X\right),Y\right] = \left[X,g\left(Y\right)\right],\label{poneta}
\end{equation}
for all $X,Y\in L$.
\end{definition}
Note that in general $C\left(L\right)\subset QC\left(L\right)$. In this section, we prove that $C\left(L\right) = QC\left(L\right)$ for $L = \mathfrak{P}$.
\begin{proposition}
Recall that $E$ is the Euler vector field and let $Y\in\mathfrak{P}_m$ with $m\geq-1$. Each element $f$ of the centroid $C\left(\mathfrak{P}\right)$ is homogeneous of degree $0.$ If we denote $f\left(E\right) = \sum_{i\geq -1} E'_i$ and $f\left(Y\right) = \sum_{i\geq -1} Y_i$, where the $E'_i,Y_i$ are in $\mathfrak{P}_i$, then $Y_{m\neq0} = \left[\frac{1}{m}E'_0, Y\right]$. Moreover $Y_0$ is completely determined by $\left[-E'_0, \left[C,Y\right]\right] = \left[C,Y_0\right] \text{for all $C \in {\mathfrak{P}_{-1}}$}$ with $\left[E'_0,Y\right] = 0$ for $m = 0$.\label{der}
\end{proposition}
\begin{proof}
If $f\in C\left(\mathfrak{P}\right)$, then we have (\ref{qs}) for $X = E$ and $Y$ an homogeneous element in $\mathfrak{P}_m$ with $m\geq-1$:
\begin{equation}
f\left[E,Y\right] = \left[fE,Y\right] = \left[E,fY\right].\label{be}
\end{equation}
Thus
\begin{equation}
mf\left(Y\right) = \sum_{i\geq -1} \left[E'_i,Y\right] = \sum_{i\geq -1}i Y_i.
\end{equation}
Then, we have $mf\left(Y\right) = \sum_{i\geq -1}m Y_i = \sum_{i\geq -1}i Y_i$. In term of homogeneous components $m Y_i = i Y_i$ and for $m\neq i$, it yields $Y_i = 0$. Hence, $f\left(Y\right) = Y_m$ and $f$ is homogeneous of degree 0. So, we get $f\left(E\right) = E'_0$ ($E$ of degree 0, then by $f$ is homogeneous of degree 0, $f\left(E\right) = E'_0$). Moreover, in term of homogeneous components of (\ref{be}), we have
\begin{equation}
\text{for all $m\neq0$},\ Y_m = \left[\frac{1}{m}E'_0,Y\right]\label{bir}
\end{equation}
In addition, we have for a constant field $C\in\mathfrak{P}_{-1}$
\begin{equation}
f\left[C,Y\right] = \left[fC,Y\right] = \left[C,fY\right].\label{bea}
\end{equation}
If $m = 0$, by the previous result: $Y$ homogeneous of degree 0, then $f(Y) = Y_0$ and $\left[C,Y\right]$ of degree -1. By (\ref{bir}) $f\left[C,Y\right] = \left[-E'_0,\left[C,Y\right]\right]$ and $fC = \left[-E'_0,C\right]$. Thus, from (\ref{bea}), we obtain
\begin{equation}
\left[-E'_0,\left[C,Y\right]\right] = \left[\left[-E'_0,C\right],Y\right] = \left[C,Y_0\right]\label{sae}.
\end{equation}
So by the Jacobi identity and (\ref{sae}),we have
\begin{equation}
\left[C,\left[-E'_0,Y\right]\right] = 0, \quad \forall C\in\mathfrak{P}_{-1}.
\end{equation}
It yields $\left[E'_0,Y\right] = 0$ for $m = 0$. Moreover,
\begin{equation}
\left[-E'_0,\left[C,Y\right]\right] = \left[C,Y_0\right] \text{for all $C\in\mathfrak{P}_{-1}$}\label{sat}
\end{equation}
determines $Y_0$, because $Y_0$ is of the form $\beta^{i}_jx^j \frac{\partial}{\partial x^i}$ with $\beta^{i}_j\in\mathbb{R}$ and $C = \frac{\partial}{\partial x^k}$ with $k = 1,\dots,n$, then $\left[C,Y_0\right] = \beta^{i}_k\frac{\partial}{\partial x^i} = \left[-E'_0,\left[C,Y\right]\right]$ determines $\beta^{i}_k$ for all $k$ and $i$.
\end{proof}
\begin{corollary}
If $\mathfrak{H}_0\subset\mathfrak{P}_{0}$ then $C\left(\mathfrak{P}\right) = \left\langle Id\right\rangle_\mathbb{R}$, where $Id$ is the identity map. If $\mathfrak{H}^d_0 = \mathfrak{P}_{0}$, $C\left(\mathfrak{P}\right) = \left\langle Id,Id_1,\dots,Id_n\right\rangle_\mathbb{R}$ where $Id_i$ with $i$ fixed is such that
\[
Id_i\left(x^i \frac{\partial}{\partial x^i}\right) = x^i \frac{\partial}{\partial x^i}, Id_i\left(\frac{\partial}{\partial x^i}\right) = \frac{\partial}{\partial x^i}
\]
and null otherwise on $\mathfrak{P}_{0}$ where $Id_i\left(Y\in\mathfrak{P}_{m\neq0}\right) = \left[\frac{1}{m}x^i \frac{\partial}{\partial x^i},Y\right]$.\label{bon}
\end{corollary}
\begin{proof}
We can refer on the previous Proposition~with its notions saying that if $f\in C\left(\mathfrak{P}\right)$, then $\left[E'_0,Y\right] = 0$ for $m = 0$. If we run $Y$ on the set of linear diagonal vector fields (that is to say $Y$ takes all values of all elements in this set), we have $E'_0 = 0$ or $E'_0 = cE+c'^{i}x^i \frac{\partial}{\partial x^i}$ with $c\neq0$ and $c'^{i}\in\mathbb{R}$. By (\ref{sat}), $\left[C,Y_0\right] = 0$ or $\left[C,Y_0-cY\right] = 0$ for all $C\in\mathfrak{P}_{-1}$ of the form $\frac{\partial}{\partial x^k}$ with $k = 1,\dots,n$ except for $\mathfrak{P}$ where $\mathfrak{P}_0 = \left\langle x^1 \frac{\partial}{\partial x^1},\dots,x^n \frac{\partial}{\partial x^n}\right\rangle_\mathbb{R}$. Then it yields $Y_0 = 0$ or $Y_0 = cY$ and (\ref{bir}) $Y_{m\neq0} = 0$ or $Y_m = cY$ for all $m\geq-1$ except for $\mathfrak{P}_0 = \left\langle x^1 \frac{\partial}{\partial x^1},\dots,x^n \frac{\partial}{\partial x^n}\right\rangle_\mathbb{R},$ otherwise $\left[E'_0,Y\right] = 0$ for all $Y\in\mathfrak{P}_0$ is not satisfied. If we have this last condition on $\mathfrak{P}_0$, we can conclude following the decomposition of $E'_0$.
\end{proof}
\begin{example}
In $\mathbb{R}^2$ with the usual coordinates system $(x,y)$, $\mathfrak{P} = \left\langle\frac{\partial}{\partial x},\frac{\partial}{\partial y},x\frac{\partial}{\partial x},y\frac{\partial}{\partial y}\right\rangle_\mathbb{R}$, $f = Id+Id_2$ in $C\left(\mathfrak{P}\right)$ is such that
\[
f\left(\frac{\partial}{\partial x}\right) = \frac{\partial}{\partial x}, \quad
f\left(\frac{\partial}{\partial y}\right) = 2\frac{\partial}{\partial y}, \quad
f\left(x\frac{\partial}{\partial x}\right) = x\frac{\partial}{\partial x}, \quad
f\left(y\frac{\partial}{\partial y}\right) = 2y\frac{\partial}{\partial y}.
\]
\end{example}
\begin{corollary}
If $\mathfrak{P}_0 = \left\langle E\right\rangle_\mathbb{R}$, then $C\left(\mathfrak{P}\right) = \left\langle Id\right\rangle_\mathbb{R}.$
\end{corollary}
\begin{proof}
If we denote $f\left(E\right) = \sum_{i\geq -1} E'_i$ with $f\in C\left(\mathfrak{P}\right)$, then $E'_0 = c E$ where $c\in\mathbb{R}$. By equations from Proposition~\ref{der}, $f\left(Y\right) = cY$ for all $Y\in\mathfrak{P}_{m\geq-1}$.
\end{proof}
\begin{corollary}
If $E'_0 = cE$ for $c\in\mathbb{R}$, then $f\in C\left(\mathfrak{P}\right)$ is always of the form $cId.$
\end{corollary}
\begin{proof}
This follows directly from the previous corollary.
\end{proof}
\begin{remark}
In general, $C\left(\mathfrak{P}\right)$ is neither $\left\langle Id\right\rangle_\mathbb{R}$ nor $\left\langle Id,Id_1,\dots,Id_n\right\rangle_\mathbb{R}$, for which we can take an example of $\mathfrak{P} = \left\langle\frac{\partial}{\partial x},\frac{\partial}{\partial y}, E, x\frac{\partial}{\partial y}\right\rangle_\mathbb{R}$ in the usual coordinates $(x,y)$ on $\mathbb{R}^2$ and $f\in C\left(\mathfrak{P}\right)$ such that $f\left(\frac{\partial}{\partial x}\right) = \frac{\partial}{\partial x}+\frac{\partial}{\partial y}, f\left(\frac{\partial}{\partial y}\right) = \frac{\partial}{\partial y},f\left(E\right) = E+ x\frac{\partial}{\partial y}, f\left(x\frac{\partial}{\partial y}\right) = x\frac{\partial}{\partial y}$.
\end{remark}
\begin{proposition}
The element $f$ of the quasicentroid $QC\left(\mathfrak{P}\right)$ is such that $f\left(E\right) = E'_0$ and $f\left(Y\right) = \left[\frac{1}{m}E'_{0},Y\right]$, where the $\left(E'_0,Y\right)$ is in $\mathfrak{H}_0\times \mathfrak{H}_{m\geq-1,m\neq0}$. Moreover, $\left[E'_{0},Y\right] = 0$ for $m = 0$, then $Y_0$ is completely defined by $\left[-E'_0,\left[C,Y\right]\right] = \left[C,Y_0\right] \text{for all $C\in\mathfrak{P}_{-1}$}$.\label{de}
\end{proposition}
\begin{proof}
Following the above notations and using the second part of the Eq. (\ref{be}), We obtain
\[
\left[f\left(E\right),E\right] = \left[E,f\left(E\right)\right], \text{so} f\left(E\right) \text{is homogeneous of degree 0}
\]
and
\[
\sum_{i\geq -1} \left[E'_i,Y\right] = \sum_{i\geq -1}i Y_i.
\]
In term of homogeneous components:
if $i\leq m-2$, $Y_i = 0$ and for $i\geq-1$,
\begin{equation}
\left[E'_i,Y\right] = (m+i) Y_{m+i}.\label{don}
\end{equation}
Because only $E'_0$ is the only term of $f\left(E\right)$, so $Y_{m}$ is the only term of $f\left(Y\right)$ and for $m\neq0$, $Y_{m} = \left[\frac{1}{m}E'_{0},Y\right]$. For $m+i = 0$ in (\ref{don}) for $m = 0$, we have $\left[E'_{0},Y\right] = 0$. Using Jacobi identity and the last equality, we get (\ref{sae}) for $f$ in the quasicentroid and it yields the last assertion of our theorem.
\end{proof}
Now, we can conclude for the centroid and the quasicentroid of $\mathfrak{P}$:
\begin{corollary}
For all $\mathfrak{P}$, $C\left(\mathfrak{P}\right) = QC\left(\mathfrak{P}\right)$.
\end{corollary}
\section{Quasiderivation of $\mathfrak{P}$}
\begin{definition}\label{sania}The set of quasiderivation $QDer\left(L\right)$ of a Lie algebra $L$ is the Lie algebra of endomorphism $f$ in $L$ such that there exists another endomorphism $g$ of $L$ such that
\begin{equation}
\left[f\left(X\right),Y\right]+\left[X,f\left(Y\right)\right] = g\left[X,Y\right]\ \forall X,Y\in L.\label{neta}
\end{equation}
\end{definition}
It is more convenient to denote $f\in QDer\left(L\right)$ by $(f,f,g)$ and $f\in QC\left(L\right)$ by $(f,-f,0)$. A quasiderivation of $\mathfrak{P}$ is a generalization of a derivation of $\mathfrak{P}$. We compute $QDer\left(\mathfrak{P}\right)$, it is clear that:
\[
QDer\left(\mathfrak{P}\right) = QDer_0\left(\mathfrak{P}\right)+\sum_{i\in \mathbb{Z}^\ast} QDer_{i}\left(\mathfrak{P}\right)
\]
where $QDer_i\left(\mathfrak{P}\right)$ is the set of quasiderivation of degree $i$.

Lemma 4.1 of~\cite{luk} can be adapted to our $\mathfrak{P}$ taking account that it doesn't need to specify the dimension of $\mathfrak{P}$ and $\mathfrak{P}$ contains a torus $\left\langle E\right\rangle_\mathbb{R}$. So
\begin{proposition}\label{gena}We have the following description of $QDer\left(\mathfrak{P}\right)$:
\[
QDer\left(\mathfrak{P}\right) = QDer_0\left(\mathfrak{P}\right)+\sum_{i\in \mathbb{Z}} QDer'_{i}\left(\mathfrak{P}\right)+ad_\mathfrak{P}
\]
where
\[
QDer'_{i}\left(\mathfrak{P}\right) = \left\{f\in QDer_{i}\left(\mathfrak{P}\right)/ f(E) = 0\right\}.
\]
\end{proposition}
In the rest of this paper, we assume $n\geq2$.
\begin{lemma}
Each element $(f,f,g)\in QDer_0\left(\mathfrak{P}\right)$ such that $f(E) = 0$, $f_{|\mathfrak{P}_{-1}} \equiv0$ and $g_{|\left[\mathfrak{P}_1,\mathfrak{P}_{-1}\right]}\equiv0$, verifies $f\equiv0$ and $g_{|\left(\mathfrak{P}_0\ominus\left[\mathfrak{P},\mathfrak{P}\right]\right)}$ is arbitrary with $g_{|\mathfrak{P}_{i\neq0}+\left[\mathfrak{P}_0,\mathfrak{P}_0\right]}\equiv0$.\label{hc}
\end{lemma}
\begin{proof}
Let us reason by induction knowing that $f$ is homogeneous of degree 0.
If $X\in\mathfrak{P}_{-1}$ then $f(X) = 0$ by hypothesis, by (\ref{neta}) where $Y = E$, we have:
\[
\left[f\left(X\right),E\right]+\left[X,f\left(E\right)\right] = g\left[X,E\right].
\]
Then $-f\left(X\right) = -g\left(X\right)$, because $f(E) = 0$ and $f\left(X\right)$ of degree -1. It yields $g(X) = 0$ from $f(X) = 0$.

For all $X\in \mathfrak{P}_{0}$ by (\ref{neta}) where $Y\in\mathfrak{P}_{-1}$:
\[
\left[f\left(X\right),Y\right]+\left[X,f\left(Y\right)\right] = g\left[X,Y\right].
\]
From the previous results on $g$ and $f$, the second term of the equality is null and then $\left[f\left(X\right),Y\right] = 0$ for all $Y\in\mathfrak{P}_{-1}$. Thus $f(X)\in\mathfrak{P}_{-1}\cap\mathfrak{P}_{0} = \{0\}$ for all $X\in \mathfrak{P}_{0}$. For $X\in\mathfrak{P}_0\cap\left[\mathfrak{P},\mathfrak{P}\right]$, we have two cases.
\begin{enumerate}
	\item
 The first is $X = \left[Y,Z\right]$ with $\left(Y,Z\right)\in\mathfrak{P}_0\times\mathfrak{P}_0$ where we deduce by definition~\ref{sania} and the previous result that $0 = \left[f\left(Y\right),Z\right]+\left[Y,f\left(Z\right)\right] = g\left[Y,Z\right] = g(X)$.
 \item
  The second case is made by the hypothesis $g\left(X\right) = 0$ where $X\in\left[\mathfrak{P}_{-1},\mathfrak{P}_1\right]$.
\end{enumerate}
	It yields $g\left(X\in\mathfrak{P}_0\ominus\left[\mathfrak{P},\mathfrak{P}\right]\right)$ is arbitrary.

For all $X\in \mathfrak{P}_{1}$ by (\ref{neta}) where $Y$ runs in $\mathfrak{P}_{-1}$:
\[
\left[f\left(X\right),Y\right]+\left[X,f\left(Y\right)\right] = g\left[X,Y\right].
\]
From $g\left(\left[\mathfrak{P}_{-1},\mathfrak{P}_1\right]\right) = 0$, the second term of the above equality is null. Because of $f(Y) = 0$, we get $\left[f\left(X\right),Y\right] = 0$ for all $Y\in\mathfrak{P}_{-1}$. Thus $f(X)\in\mathfrak{P}_{-1}\cap\mathfrak{P}_{1} = \{0\}$ for all $X\in \mathfrak{P}_{1}$. But (\ref{neta}) where $Y = E$ says $f(X) = g(X)$, then $g(X) = 0$ also.

We set $k\geq2$, we suppose that $f_{|\mathfrak{P}_{t}}\equiv g_{|\mathfrak{P}_{t}}\equiv 0$ for all $t\geq k-1$. For $X\in \mathfrak{P}_{k\geq2}$, we write (\ref{neta}) with $Y\in\mathfrak{P}_{-1}$. Then we obtain $f(X)\in\mathfrak{P}_{-1}\cap\mathfrak{P}_{k\geq2} = \{0\}$ for all $X\in \mathfrak{P}_{k}$ knowing that $f_{|\mathfrak{P}_{-1}}\equiv g_{|\mathfrak{P}_{k-1}}\equiv 0$. Again by (\ref{neta}) with $Y = E$, we have $0 = -kf(X) = -kg(X) = 0$ which yields $g(X) = 0$ by $f(X) = 0$.
\end{proof}
\begin{proposition}\label{asif}Suppose that all elements of $\left[\mathfrak{P}_1,\mathfrak{P}_{-1}\right]$ are also elements of $\left[\mathfrak{P}_{0},\mathfrak{P}_0\right]$. The $(f,f,g)\in QDer_0\left(\mathfrak{P}\right)$ vanishing on $E$ is of the following form:
$\left(L_X,L_X,L_X+k\right)$, where $X\in \mathfrak{H}_0$ and $k\in End\left(\mathfrak{P},\mathfrak{P}\right)$ such that $k$ is arbitrary on $\mathfrak{P}_0\ominus\left[\mathfrak{P},\mathfrak{P}\right]$ and null elsewhere.\label{any}
\end{proposition}
\begin{proof}
Let $(f,f,g)\in QDer_0\left(\mathfrak{P}\right)$, such that $f(E) = 0$. In the system of coordinates on $\mathbb{R}^n$ $(x^1,\dots,x^n)$, we write $f\left(X_0\right) = X_1\in\mathfrak{P}_{-1}$ for $X_0\in\mathfrak{P}_{-1}$, with $f(E) = 0$. Like in the proof of Proposition~2.8 of~\cite{sam}, there exists an unique $X\in \mathfrak{H}_0$ such that $f\underset{\left\langle E\right\rangle_\mathbb{R}+\mathfrak{P}_{-1}}{ = }L_X$. Then we write the following $(f,f,g)-\left(L_X,L_X,L_X\right)\in QDer_0\left(\mathfrak{P}\right)$ with $f(E) = L_X\left(E\right) = 0$. If we denote $f' = f-L_X$ and $g' = g-L_X$, we have $f'(E) = 0$, $f'\underset{\mathfrak{P}_{-1}}{ = }0$ and $g'_{|\left[\mathfrak{P}_1,\mathfrak{P}_{-1}\right]}\equiv0$ by the first hypothesis in the statement of the present proposition. That is to say, if $X_2\in\left[\mathfrak{P}_1,\mathfrak{P}_{-1}\right]$, then there are $(Y,Z)\in\mathfrak{P}_{0}$ such that $X_2 = \left[Y,Z\right]$. Then
\[
g'\left(X_2\right) = \left[f'\left(Y\right),Z\right]+\left[Y,f'\left(Z\right)\right].
\]
By the proof of Lemma~\ref{hc}, $f'\left(Y\right) = f'\left(Z\right) = 0$, thus $g'\left(X_2\right) = 0$. By the same lemma, $(f',f',g')\equiv0$ except in $\mathfrak{P}_0\ominus\left[\mathfrak{P},\mathfrak{P}\right]$ where $g'$ is arbitrary.
\end{proof}
\begin{proposition}
Let $(f,f,g)\in QDer'_{i\geq1}\left(\mathfrak{P}\right)$ and adopt all hypothesis of Lemma~\ref{hc}. Then $f\equiv g_{|\mathfrak{P}_{j\neq0}+\left[\mathfrak{P}_0,\mathfrak{P}_0\right]}\equiv0$ and $g_{|\left(\mathfrak{P}_0\ominus\left[\mathfrak{P},\mathfrak{P}\right]\right)}$ is arbitrary. If $(f,f,g)\in QDer'_{-1}\left(\mathfrak{P}\right)$, we suppose all hypothesis of Lemma~\ref{hc} except for $f_{|\mathfrak{P}_{-1}}\equiv0$, we have the same result as the previous one.\label{bouzo}
\end{proposition}
\begin{proof}
We prove this Proposition~by induction.
Let $(f,f,g)\in QDer'_{i\geq1}\left(\mathfrak{P}\right)$.
By (\ref{neta}) where $X = E$ and $Y$ runs in $\mathfrak{P}_{-1}$, we can say that $g_{|\mathfrak{P}_{-1}}\equiv 0$ knowing $f_{|\mathfrak{P}_{-1}}\equiv 0$ and $f(E) = 0$.
If $X\in\mathfrak{P}_{0}$, (\ref{neta}) with $Y = E$ where $f(E) = 0$ gives $\left[fX,E\right] = 0$. It yields $f(X) = 0$ because $f(X)\in \mathfrak{P}_{i\neq0}$. Then, in the turn of $g$, if $X\in\left[\mathfrak{P}_{0},\mathfrak{P}_{0}\right]$, then Eq. (\ref{neta}) yields $g(X) = 0$. In addition, we have $g_{|\left[\mathfrak{P}_{1},\mathfrak{P}_{-1}\right]}\equiv 0$ by hypothesis. Thus if $X\in\left[\mathfrak{P}_{1},\mathfrak{P}_{-1}\right]$, $g(X) = 0$.
Next, $X\in\mathfrak{P}_{1}$, if $Y$ runs in $\mathfrak{P}_{-1}$, in (\ref{neta}) we obtain $\left[f\left(X\right),\mathfrak{P}_{-1}\right] = \{0\}$. But the degree of $f$ is not -2, then $f\left(X\right) = 0$. It follows that the equation (\ref{neta}) with $Y = E$ and this $X$ permit us to say that $g(X) = 0$.
Now, we suppose that $f = g = 0$ on $\mathfrak{P}_{t\leq l}$ with $l\geq1$ and $t\geq1$. Consider $X\in \mathfrak{P}_{l+1}$, if we run $Y\in\mathfrak{P}_{-1}$ on (\ref{neta}), we have $\left[f\left(X\right),\mathfrak{P}_{-1}\right] = \{0\}$. It yields $f(X) = 0$ because $f(X)$ is of degree $i+l+1$ different to $-1$. If we write $Y = E$ with this $X$ on (\ref{neta}), we have $-(l+1)g(X) = \left[f(X),E\right] = -(i+l+1)f(X)$ which is equal to 0 by the above result on $f$. Then, we get $g(X) = 0$ because $l+1\neq0$.
For $(f,f,g)\in QDer'_{-1}\left(\mathfrak{P}\right)$, it is obvious that $f_{|\mathfrak{P}_{-1}} = g_{|\mathfrak{P}_{-1}}\equiv 0.$ The rest of the proof is the same as for $(f,f,g)\in QDer'_{i\geq1}\left(\mathfrak{P}\right)$.
\end{proof}
\begin{proposition}
Consider all hypothesis of Proposition~\ref{asif}. Every $(f,f,g)$ in the set $QDer_{i\geq-1,\neq0,1}\left(\mathfrak{P}\right)$ is of the form
\[
\left(L_{\frac{-1}{i}{f(E)}_{i}},L_{\frac{-1}{i}{f(E)}_{i}},L_{\frac{-1}{i}{f(E)}_{i}}+k\right),
\]
where ${f(E)}_{i}$ the homogeneous term $f(E)$ of degree $i$ and $k$ is an endomorphism of $\mathfrak{P}$, arbitrary on $\mathfrak{P}_0\ominus\left[\mathfrak{P},\mathfrak{P}\right]$ and null elsewhere. \label{uni}
\end{proposition}
\begin{proof}
Let $(f,f,g)\in QDer_{i\geq-1,\neq0,1}\left(\mathfrak{P}\right)$ be such quasiderivation, we have
\[
\left(f' = f-L_{\frac{-1}{i}{f(E)}_{i}},f' = f-L_{\frac{-1}{i}{f(E)}_{i}},g' = g-L_{\frac{-1}{i}{f(E)}_{i}}\right)\in QDer'_{i\geq-1,\neq0,1}\left(\mathfrak{P}\right).
\]
Moreover, by using $X = E$ and $Y\in\mathfrak{P}_0$ in Eq. (\ref{neta}), we have
\[
g'\left[E,Y\right] = \left[E,f'(Y)\right].
\]
This yields $0 = \left[E,f'(Y)\right]$ and $f'(Y) = 0$ because the degree of $f'(Y)$ is nonzero. Then $g'_{|\left[\mathfrak{P}_{0},\mathfrak{P}_0\right]} = \{0\}$ from (\ref{neta}). Knowing that all elements of $\left[\mathfrak{P}_1,\mathfrak{P}_{-1}\right]$ are elements of $\left[\mathfrak{P}_{0},\mathfrak{P}_0\right]$, we have $g'_{\left[\mathfrak{P}_1,\mathfrak{P}_{-1}\right]}\equiv0$. Next, we will check whether $f'_{\mathfrak{P}_{-1}}\equiv 0$. If $i = -1$, it is obvious that $f'_{\mathfrak{P}_{-1}} = g'_{\mathfrak{P}_{-1}}\equiv0$. If $i\geq2$, we take $X = \frac{\partial}{\partial x^t}\in\mathfrak{P}_{-1}$ where $1\leq t\leq n$ is fixed. We write $f'\left(X\right) = \underset{k\neq t}{P^k}\frac{\partial}{\partial x^k}+P^t \frac{\partial}{\partial x^t}$ where all $P^s$ is a polynomial of degree $i$, in which all $P^s$ doesn't depend on $x^s$ because of the hypothesis that each linear diagonal field doesn't belong to $\left[\mathfrak{P},\mathfrak{P}\right]$. When $Y = x^t \frac{\partial}{\partial x^t}$ in Eq. (\ref{neta}), $0 = g'\left[\frac{\partial}{\partial x^t}, x^l \frac{\partial}{\partial x^l}\right] = \left[f'(X), x^l \frac{\partial}{\partial x^l}\right]$ with a fixed $l\neq t$. Then
\[
P^l \frac{\partial}{\partial x^l}+\beta^k \underset{k\neq t,l}{P^k}\frac{\partial}{\partial x^k}+\alpha^t P^t \frac{\partial}{\partial x^t} = 0,
\]
where $\alpha^t, \beta^k \in\mathbb{R}$. It yields $P^l = 0$ for all $l\neq t$. Thus,
\[
g'\left(\frac{\partial}{\partial x^t}\right) = g'\left[\frac{\partial}{\partial x^t}, x^t \frac{\partial}{\partial x^t}\right] = \left[f'\left(\frac{\partial}{\partial x^t}\right),x^t \frac{\partial}{\partial x^t}\right] = f'\left(\frac{\partial}{\partial x^t}\right).
\]
When in Eq.(\ref{neta}), $X\in\mathfrak{P}_{l\geq-1}$ and $Y = E$, we get
\[
\frac{-l}{i}g(X) = L_{\frac{-1}{i}{f(E)}_{i}}(X)-\frac{i+l}{i}f(X).
\]
Then $\left(f-L_{\frac{-1}{i}{f(E)}_{i}}\right)(X) = \frac{f\left(X\right)-g\left(X\right)}{i}$ where $X\in\mathfrak{P}_{-1}$. If $i\geq2$, $f'\left(X\right)-g'\left(X\right) = 0$ for $X\in\mathfrak{P}_{-1}$ by the above result, then $f\left(X\right)-g\left(X\right) = 0$.
Then by Proposition~\ref{bouzo}, we have $f'\equiv g'_{|\mathfrak{P}_{j\neq0}+\left(\mathfrak{P}_0\cap\left[\mathfrak{P},\mathfrak{P}\right]\right)}\equiv0$ and $g'_{|\left(\mathfrak{P}_0\ominus\left[\mathfrak{P},\mathfrak{P}\right]\right)}$ is arbitrary.
\end{proof}
\begin{proposition}
We suppose all hypothesis of Proposition~\ref{any}. Every $(f,f,g)$ in the set $QDer_{1}\left(\mathfrak{P}\right)$ is of the form
\[
\left(-L_{{f(E)}_{1}},-L_{{f(E)}_{1}},-L_{{f(E)}_{1}}+k\right)+\left(f'',f'',0\right)
\]
where ${f(E)}_{1}$ is the homogeneous term of $f(E)$ of degree $1$, $k$ is arbitrary on $\mathfrak{P}_0\ominus\left[\mathfrak{P},\mathfrak{P}\right]$ and null elsewhere, $f''_{|\mathfrak{P}_{l\neq-1}}$ is null and $f''\left(\mathfrak{P}_{-1}\right)$ is a subset of $\left\langle x^1 \frac{\partial}{\partial x^1},\dots, x^n \frac{\partial}{\partial x^n}\right\rangle_\mathbb{R}$. \label{unie}
\end{proposition}
\begin{proof}
Let $(f,f,g)\in QDer_{1}\left(\mathfrak{P}\right)$ be such quasiderivation, we write for $i = 1$:
\[
\left(f' = f-L_{\frac{-1}{i}{f(E)}_{i}},f' = f-L_{\frac{-1}{i}{f(E)}_{i}},g' = g-L_{\frac{-1}{i}{f(E)}_{i}}\right)\in QDer'_{i\neq0}\left(\mathfrak{P}\right).
\]
We can assume all results of the first part of the previous proof ending on $g'_{|\left[\mathfrak{P}_1,\mathfrak{P}_{-1}\right]}\equiv0$ with $i = 1$. In this result, we have $f'_{|\mathfrak{P}_0}\equiv0$.

Let us take $X\in\mathfrak{P}_{-1}$. When $Y = E$ in (\ref{neta}), $g'(X) = \left[f'(X),E\right] = 0$ because $f'(X)$ of degree 0 and $f'(E) = 0$.
Let us precise $X = \frac{\partial}{\partial x^j}$, $\left[\frac{\partial}{\partial x^j},x^j \frac{\partial}{\partial x^j}\right] = X$ and by Eq. (\ref{neta}) with $Y = x^j \frac{\partial}{\partial x^j}$, $g'\left(X\right) = \left[f'(X),x^j \frac{\partial}{\partial x^j}\right] = 0$ knowing that $f'(x^j \frac{\partial}{\partial x^j}) = 0$. For all $k\neq j$ fixed, $\left[\frac{\partial}{\partial x^j},x^k \frac{\partial}{\partial x^k}\right] = 0$, then (\ref{neta}) where $Y = x^k \frac{\partial}{\partial x^k}$ yields $\left[f\left(\frac{\partial}{\partial x^j}\right),x^k \frac{\partial}{\partial x^k}\right] = 0$ for all $k\neq j$. We can conclude that $f'\left(\frac{\partial}{\partial x^j}\right) = \alpha^t x^t \frac{\partial}{\partial x^t}$. If we do $g'\left[\frac{\partial}{\partial x^j},\frac{\partial}{\partial x^l}\right] = \left[\alpha^t x^t \frac{\partial}{\partial x^t},\frac{\partial}{\partial x^l}\right]+\left[\frac{\partial}{\partial x^j},\beta^s x^s \frac{\partial}{\partial x^t}\right]$, we obtain $\alpha^l = \beta^j = 0$ for all $j\neq l$. Then $f'\left(\frac{\partial}{\partial x^j}\right) = \alpha^j x^j \frac{\partial}{\partial x^j}$. Thus, let $f''\in End\left(\mathfrak{P}\right)$ be such that $f''$ is null except on $\mathfrak{P}_{-1}$ where the value is on $\left\langle x^1 \frac{\partial}{\partial x^1},\dots, x^n \frac{\partial}{\partial x^n}\right\rangle_\mathbb{R}$ such that $f'\left(\frac{\partial}{\partial x^j}\right) = f''\left(\frac{\partial}{\partial x^j}\right)$. It is easy to verify that $(f'',f'',0)$ is in $QDer'_{1}\left(\mathfrak{P}\right)$.

We get
\[
\left(f''' = f'-f'',f''' = f'-f'',g'\right)\in QDer'_{1}\left(\mathfrak{P}\right)
\]
verifying $f'''(E) = 0$, $f'''_{|\mathfrak{P}_{-1}}\equiv0$ and $g'_{|\left[\mathfrak{P}_1,\mathfrak{P}_{-1}\right]}\equiv0$. Then we can use Proposition~\ref{bouzo} and have the following results: $f''' = 0$, $g'_{|\mathfrak{P}_{t\neq0}+\left(\mathfrak{P}_0\cap\left[\mathfrak{P},\mathfrak{P}\right]\right)}\equiv0$ and $g'_{|\left(\mathfrak{P}_0\ominus\left[\mathfrak{P},\mathfrak{P}\right]\right)}$ is arbitrary. It follows the final results.
\end{proof}
\begin{definition}
For a fixed $i_0$, $j_1< j_2 \dotsb < j_p$ with $1\leq p\leq n$, $\left(\alpha_{j_1},\dots,\alpha_{j_p}\right)\in\mathbb{N}^p$ and $j_l\in\{1,\dots,n\}$, the vector field ${\left(x^{j_1}\right)}^{\alpha_{j_1}}\dots {\left(x^{j_p}\right)}^{\alpha_{j_p}}\frac{\partial}{\partial x^{i_0}}$ is called monomial. We recall that $ \mathfrak {P} $ is separated, if all homogeneous elements of degree $k\geq-1$ are spanned by monomial vector field in $ \mathfrak {P} $.
\end{definition}
\begin{remark}
In the proof of the previous theorem, the following fact is true, $x^j \frac{\partial}{\partial x^j}\notin\left[\mathfrak{P},\mathfrak{P}\right]$ if we adopt the hypothesis on $\mathfrak{P}$: $\left[\mathfrak{P}_1,\mathfrak{P}_{-1}\right]\subset\left[\mathfrak{P}_0,\mathfrak{P}_{0}\right]$. In the same proof, we take $X = \frac{\partial}{\partial x^j}$. It is supposed that $\mathfrak{P}$ is separated. We can reason as the following, if there exists $Y = x^k \frac{\partial}{\partial x^j}\in\mathfrak{P}_0$ or $Y = x^j \frac{\partial}{\partial x^k}\in\mathfrak{P}_0$ for $j\neq k$, $0 = g'\left[X,Y\right] = \left[f'\left(X\right),Y\right]$ from (\ref{neta}) knowing $f'(Y) = 0$ and $g'_{|\mathfrak{P}_{-1}}\equiv0$. It yields $\alpha^j = 0$.
\end{remark}
\begin{example}
We are in $\mathbb{R}^2$ with the habitual coordinates $(x,y)$, $\mathfrak{P} = \left\langle\frac{\partial}{\partial x},\frac{\partial}{\partial y},x\frac{\partial}{\partial x},y\frac{\partial}{\partial y} \right\rangle_\mathbb{R}$. Define an endomorphism $f''$ on $\mathfrak{P}$, such that $f''\left(\frac{\partial}{\partial x}\right) = x\frac{\partial}{\partial x}$, $f''\left(\frac{\partial}{\partial y}\right) = y\frac{\partial}{\partial y}$ and null otherwise. Then $(f'',f'',0)$ is an example of quasiderivation of type Proposition~\ref{unie}, which is neither a derivation nor an element of $QC\left(\mathfrak{P}\right)$.\label{fs}
\end{example}
\begin{proposition}
Let $(f,f,g)\in QDer'_{i\leq-2}\left(\mathfrak{P}\right)$, then $f_{|\mathfrak{P}_{t\neq-(i+1)}}\equiv g_{|\mathfrak{P}_{t\neq-(i+1)}}\equiv0$ with $f(X) = (i+1)g(X)$ for all $X\in \mathfrak{P}_{-(i+1)}$ and $g_{|\left(\mathfrak{P}_0\ominus\left[\mathfrak{P},\mathfrak{P}\right]\right)}$ is arbitrary.\label{bouz}
\end{proposition}
\begin{proof}
Let $(f,f,g)\in QDer'_{i\leq-2}\left(\mathfrak{P}\right)$ a such quasiderivation. By the fact that $\mathfrak{P}_{t\leq-2} = \{0\}$, $f = g\equiv0$ on $\mathfrak{P}_{t\leq0}$. We reason by induction, we suppose that $f = g = 0$ on $\mathfrak{P}_{t\leq 0}$. Consider $X\in \mathfrak{P}_{t+1}$ with $t\geq0$, if we write $Y = E$ on (\ref{neta}), we have
\[
-(t+1)g(X) = \left[f(X),E\right] = -(i+t+1)f(X).
\]
We take (\ref{neta}) and run $Y$ on $\mathfrak{P}_{-1}$, we obtain $\left[f\left(X\right),Y\right] = \{0\}$. Hence, if $i+t\neq-2$, then $f\left(X\right) = 0$ and $g(X) = 0$ by the previous equation. Otherwise, $f(X) = -(t+1)g(X)$ for all $X\in \mathfrak{P}_{t+1}$, that is to say in the case where $t+1 = -1-i$.
\end{proof}
\begin{proposition}
If each element of $\left[\mathfrak{P}_1,\mathfrak{P}_{-1}\right]$ is an element of $\left[\mathfrak{P}_{0},\mathfrak{P}_0\right]$ and $\mathfrak{P}$ is separated with all linear diagonal fields in $\mathfrak{P}$, then $(f,f,g)\in QDer'_{i\leq-2}\left(\mathfrak{P}\right)$ is null except $g_{|\left(\mathfrak{P}_0\ominus\left[\mathfrak{P},\mathfrak{P}\right]\right)}$ which is arbitrary.\label{ita}
\end{proposition}
\begin{proof}
To fix our idea, we cannot take $X = x^1 P(x^{t\neq1})\frac{\partial}{\partial x^1}\in\mathfrak{P}_{l\geq1}$ because in this case $x^1 \frac{\partial}{\partial x^1}\in\mathfrak{P}_0\cap\left[\mathfrak{P}_1,\mathfrak{P}_{-1}\right] \ \setminus \ \left[\mathfrak{P}_{0},\mathfrak{P}_0\right]$. So we are forced to take $X = P(x^{t\neq1})\frac{\partial}{\partial x^1}\in\mathfrak{P}_{l\geq1}$. By using $f(X) = (i+1)g(X)$ for $X\in \mathfrak{P}_{l = -(i+1)}$ where $-(i+1)\geq1$ from Proposition~\ref{bouz} and when $Y = x^1 \frac{\partial}{\partial x^1}$ in the Eq. (\ref{neta}) with $g(X) = \alpha^j \frac{\partial}{\partial x^j}$ where $\alpha^j \in\mathbb{R}$: we obtain
\[
\alpha^j \frac{\partial}{\partial x^j} = (i+1)\alpha^1 \frac{\partial}{\partial x^1}.
\]
It yields $\alpha^j = 0$ for all $j\neq1$ and $\alpha^1 = 0$ because $i+1\leq-1$. Thus $g(X)$ and $f(X)$ are zero.
\end{proof}
Next, by using result on odd-derivations in~\cite{pri}, we can prove that $(f, f, g)\in QDer'_{-2}\left(\mathfrak{P}\right)$ is an odd-derivations of $\mathfrak{P}$ in the case $\mathfrak{P}_0 = \left\langle E\right\rangle_\mathbb{R}\subsetneq\left[\mathfrak{P}_{1},\mathfrak{P}_{-1}\right]$.
\begin{proposition}
Let $(f,f,g)\in QDer'_{-2}\left(\mathfrak{P}\right)$. Then $g_{|\mathfrak{P}_0\cap\left[\mathfrak{P}, \mathfrak{P}\right]}\equiv 0$ and $g_{|\left(\mathfrak{P}_0\ominus\left[\mathfrak{P},\mathfrak{P}\right]\right)}$ is arbitrary. In addition, $g_{|\mathfrak{P}_{i\neq0,1}}\equiv0$ with $g\left(X\right) = -f\left(X\right)$ for $X\in\mathfrak{P}_1$, $f\left[Y,X'\right] = \left[fX',Y\right]$ for all $\left(X',Y\right)\in\mathfrak{P}_{2}\times\mathfrak{P}_{-1}$ and $f_{|\mathfrak{P}_{i\neq1,2}}\equiv0$.\label{bol}
\end{proposition}
\begin{proof}
If $f$ have degree $-2$, then $f_{|\mathfrak{P}_{i\leq0}}\equiv0$. Thus, applying (\ref{neta}) with $(X, Y)\in\mathfrak{P}_{1}\times\mathfrak{P}_{-1}$ and $(X, Y)\in\mathfrak{P}_{0}\times\mathfrak{P}_{0}$ resp., knowing that $f$ is of degree $-2$, we obtain that $g\left[X,Y\right] = 0$ but $g_{|\left(\mathfrak{P}_0\ominus\left[\mathfrak{P},\mathfrak{P}\right]\right)}$ is arbitrary. Now, let's take $X\in\mathfrak{P}_{1}$ and $Y = E$ in (\ref{neta}), we deduce that $g\left(X\right) = -f\left(X\right)$ for $f$ of degree $-2$. Next, we consider $X\in\mathfrak{P}_{2}$ and $Y = E$ in Eq. (\ref{neta}), we find $g(X) = 0$ with $f(X)$ is determined by the same relation, if $X\in\mathfrak{P}_{2}$ and $Y\in\mathfrak{P}_{-1}$, we find that $f\left[Y, X\right] = \left[f(X), Y\right]$ taking account the value of $f$ and $g$ on $\mathfrak{P}_{1}$. If $X\in\mathfrak{P}_{k\geq3}$ and $Y$ runs in $\mathfrak{P}_{-1}$, by induction we have $f\left(X\right) = 0$ and $g\left(X\right) = 0$ by Eq. (\ref{neta}) when $Y = E$.
\end{proof}
Let us recall the Theorem 3.10 of~\cite{pri}:
\begin{theorem}\label{mao}We suppose that $m\geq3$ is odd. A linear map $D$ of degree $-2$ on $\mathfrak{P}$, null on $\mathfrak{P}\ominus\mathfrak{P}_1$ is an $m$-derivation, if and only if:
\begin{equation}
\left[D\left(\mathfrak{P}_1\right),\mathfrak{P}-\mathfrak{P}_0\right] = \{0\};\label{no}
\end{equation}
\begin{equation}
D\left[X,Y\right] = \left[D\left(Y\right),X\right]\ \forall\left(X,Y\right)\in\mathfrak{P}_0\times\mathfrak{P}_1\label{ak};
\end{equation}
\begin{equation}
\left[D\left(\mathfrak{P}_1\right),\left[\mathfrak{P},\mathfrak{P}\right]\cap\mathfrak{P}_0\right] = \{0\}\label{ka};
\end{equation}
if $\left[X_1,\left[X_2,\dots,\left[X_i,\dots,\left[X_{m-1},X_m\right]\dots\right]\right]\right]\in\mathfrak{P}_{1}$ where $i$ is the first index with $X_i\in\mathfrak{P}_{1}$ with the existence of $1\leq j<i$ such that $X_j\in\mathfrak{P}_{-1}\cup\mathfrak{P}_{t\geq2}$, then
\begin{equation}
D\left[X_1,\left[X_2,\dots,\left[X_i,\dots,\left[X_{m-1},X_m\right]\dots\right]\right]\right] = 0\label{aak}
\end{equation}
\end{theorem}
\begin{proposition}
If $(f,f,g)\in QDer'_{-2}\left(\mathfrak{P}\right)$ and if for all $X\in\mathfrak{P}_2$ there exists $X'\in\mathfrak{P}_{-1}$ such that $\left[X,X'\right] = 0$. Moreover, suppose that $\mathfrak{P}$ is separated and $\mathfrak{P}_0 = \left\langle E\right\rangle_\mathbb{R}\subsetneq\left[\mathfrak{P}_{1},\mathfrak{P}_{-1}\right]$. Then the endomorphism $f$ is an odd-derivation of $\mathfrak{P}$, the endomorphism $g = -f$ on $\mathfrak{P}_{1}$ and vanishes otherwise except on $\left\langle E\right\rangle_\mathbb{R}$ where $g$ is arbitrary.\label{pre}
\end{proposition}
\begin{proof}
The $f$ is an $\mathbb{R}$-linear map of degree $-2$ on $\mathfrak{P}$ with $(f,f,g)\in QDer'_{-2}\left(\mathfrak{P}\right)$. For all $X'\in\mathfrak{P}_{2}$, $\left[X',E\right] = -2X'$, then $\left[f\left(X'\right), E\right]+\left[X', f(E)\right] = g\left[X', E\right] = -2g\left(X'\right) = 0$ by (\ref{neta}) and Proposition~\ref{bol} saying $f(E) = 0$ and $g(X') = 0$. It yields that there exists $c\in\mathbb{R}$ such that $f\left(X'\right) = cE$. In addition, for all $\left(X', Y\right)\in\mathfrak{P}_{2}\times\mathfrak{P}_{-1}$, we have $f\left[Y, X'\right] = \left[fX',Y\right]$ using Proposition~\ref{bol}. By hypothesis, there is a $Y\in\mathfrak{P}_{-1}$ such that $\left[Y, X'\right] = 0$. We can choose $Y = \frac{\partial}{\partial x^j}\in\mathfrak{P}_{-1}$ for a fixed $j$. It conducts to $0 = \left[fX',Y\right] = -c\frac{\partial}{\partial x^j}$ and $c = 0$. We can conclude that $f_{|\mathfrak{P}_{2}}\equiv0$. Now, if we look at Theorem~\ref{mao}, we must check the $4$ relations in that theorem for $f$. Let us take $X\in\mathfrak{P}_{1},\ Y\in\mathfrak{P}_{i\neq0},\ Z\in\mathfrak{P}_{0}$:
\begin{enumerate}
	\item $\left[f\left(X\right),Y\right] = g\left[X,Y\right]-\left[X,f(Y)\right]$ by definition of quasiderivation. Since we have $\left[X,Y\right]\notin\mathfrak{P}_{1}$, then $\left[f\left(X\right),Y\right] = -\left[X,f(Y)\right]$ by Proposition~\ref{bol} saying $g_{|\mathfrak{P}_0\cap\left[\mathfrak{P}, \mathfrak{P}\right]}\equiv 0$ and $g_{|\mathfrak{P}_{i\neq0,1}}\equiv0$. Then if $Y\notin\mathfrak{P}_{1}$, $\left[f\left(X\right),Y\right] = -\left[X,f(Y)\right] = 0$ by Proposition~\ref{bol} and the above result on $f_{|\mathfrak{P}_{2}}\equiv0$. But if $Y\in\mathfrak{P}_{1}$, because of $\left[f\left(X\right),Y\right]$ and $\left[X,f(Y)\right]$ are in $\mathfrak{P}_0 = \left\langle E\right\rangle_\mathbb{R}\subsetneq\left[\mathfrak{P}_{1},\mathfrak{P}_{-1}\right]$, we have $\left[f\left(X\right),Y\right] = \left[X,f(Y)\right] = 0.$ So we can resume $\left[f\left(X\right),Y\right] = 0$ for all $Y\in\mathfrak{P}_{i\neq0}$.
	\item We have $g\left[Z,X\right] = \left[Z,f\left(X\right)\right]$, taking into account that $f = -g$ on $\mathfrak{P}_{1}$, we obtain $f\left[Z,X\right] = \left[f\left(X\right),Z\right]$.
	\item
\begin{itemize}
\item(i) First, consider $X',X''\in\mathfrak{P}_{0}$, $\left[f\left(X\right),\left[X',X''\right]\right] = g\left[X,\left[X',X''\right]\right].$ Then it is equal to $g\left[\left[X,X'\right],X''\right]+g\left[X',\left[X,X''\right]\right]$ by Jacobi identity. The property of $f$ conducts this value to $\left[f\left[X,X'\right],X''\right]+\left[X',f\left[X,X''\right]\right]$ and by $f = -g$ in $\mathfrak{P}_{1}$, it is of the following form
\[
\left[-g\left[X,X'\right],X''\right]+\left[X',-g\left[X,X''\right]\right].
\]
Then this expression becomes $-\left[\left[fX,X'\right],X''\right]-\left[X',\left[fX,X''\right]\right]$. By the identity of Jacobi, it is
	
\[
-\left[\left[fX,X''\right],X'\right]-\left[\left[X',fX\right],X''\right]-2\left[fX,\left[X',X''\right]\right]
\]
	and then
\[
-\left[g\left[X,X''\right],X'\right]-\left[g\left[X',X\right],X''\right]-2g\left[X,\left[X',X''\right]\right].
\]
	By $f = -g,$ we have
\[
g\left[\left[X,X''\right],X'\right]+g\left[\left[X',X\right],X''\right]-2g\left[X,\left[X',X''\right]\right]
\]
and by the Jacobi identity, we have
\[
-g\left[X,\left[X',X''\right]\right]-2g\left[X,\left[X',X''\right]\right].
\]
By the first relation, it is equal to $g\left[X,\left[X',X''\right]\right]$, thus $g\left[X,\left[X',X''\right]\right] = 0$ and $\left[f\left(X\right),\left[X',X''\right]\right] = 0.$
	\item(ii) Second, consider $\left(X',X''\right)\in\mathfrak{P}_{1}\times\mathfrak{P}_{-1}$,
	
\[
\left[f\left(X\right),\left[X',X''\right]\right] = g\left[X,\left[X',X''\right]\right].
\]
	And it becomes
\begin{equation*}
\begin{aligned}g\left[\left[X,X'\right],X''\right]+g\left[X',\left[X,X''\right]\right] &= \left[f\left[X,X'\right],X''\right]+\left[fX',\left[X,X''\right]\right]\\
&= \left[\left[fX',X\right],X''\right]+\left[X,\left[fX',X''\right]\right].
\end{aligned}
\end{equation*} This right term of the above equality becomes
	$\left[\left[fX',X\right],X''\right]$ because the degree of $f$ is $-2$. Thus $\left[\left[fX',X\right],X''\right] = 0$ by (1).
	So we conclude that
\[
\left[f\left(\mathfrak{P}_{1}\right),\left[\mathfrak{P},\mathfrak{P}\right]\cap\mathfrak{P}_{0}\right] = \{0\}.
\]
\end{itemize}
	\item Consider $2k+1$ elements $X_1, \dots, X_j,\dots, X_i, \dots, X_{2k+1}$ on $\mathfrak{P}$ as in Theorem~\ref{mao}, and knowing that $f = -g$ on $\mathfrak{P}_{1}$, we have
	
\begin{align*}
	&f\left[X_1,\left[X_2,\dots,\left[X_i,\dots,\left[X_{m-1},X_m\right]\dots\right]\right]\right]\\
&= \left[X_1,\left[X_2,\dots,-f\left[X_j,\dots,\left[X_i,\dots,\left[X_{m-1},X_m\right]\dots\right]\right]\right]\right].
\end{align*}
With the property of $g$, it is
	
\[
\left[X_1,\left[X_2,\dots,\left[X_j,f\left[X_{j+1},\dots,\left[X_i,\dots,\left[X_{m-1},X_m\right]\dots\right]\right]\right]\right]\right].
\]
By the (1), this expression vanishes.
\end{enumerate}
We conclude that $f$ is a $(2k+1)$-derivation on $\mathfrak{P}$ like in Theorem~\ref{mao}. We have $f = -g$ except on $\left\langle E\right\rangle_\mathbb{R}$ where $g$ is arbitrary.
\end{proof}
\begin{remark}
In general, the $(2k+1)$-derivation of the Proposition~\ref{pre} is not null. We can see it with one example from~\cite{pri}:
In $\mathbb{R}^3$, the Lie algebra $\mathfrak{P}$ spanned by $ \frac{\partial} {\partial x}$, $\frac{\partial} {\partial y}$, $\frac{\partial} {\partial z}$,$E$, $x \frac{\partial} {\partial z}$ and $(x) ^ 2 \frac{\partial} {\partial z}$ admits a $3$-derivation of $\mathfrak{P}$, $D_0$ defined by $D_0\left((x) ^ 2 \frac{\partial} {\partial z}\right) = \frac{\partial} {\partial z}$ and null otherwise. It is easy to check that $D_0$ is in $QDer'_{-2}\left(\mathfrak{P}\right)$. That is to say $(D_0,D_0,-D_0)\in QDer'_{-2}\left(\mathfrak{P}\right)$.
\end{remark}
\begin{proposition}
If we have the hypothesis of Proposition~\ref{any} and if all linear diagonal fields are in $\mathfrak{P}$, then every elements of $(f,f,g)\in QDer_0\left(\mathfrak{P}\right)$ is a sum of an element of $\left\langle Id,Id_1,\dots,Id_n\right\rangle_\mathbb{R}$ and a quasiderivation of the type Proposition~\ref{any}.\label{en}
\end{proposition}
\begin{proof}
If we take such $f$, we can apply Eq.(\ref{neta}) with $X = E$ and $Y$ a linear diagonal fields and find that $\left[f\left(E\right),Y\right] = 0$ for all $Y$. It conducts to the following facts, there exists $Y'$ a non-null linear diagonal field such that $f(E) = Y'$. We take $h\in \left\langle Id,Id_1,\dots,Id_n\right\rangle_\mathbb{R}$ with $\left\langle Id,Id_1,\dots,Id_n\right\rangle_\mathbb{R}\supset QC\left(\mathfrak{P}\right)$ by corollary~\ref{bon} such that $h\left(E\right) = Y'$. This implies that $(f-h,f+h,g)\in QDer'_{0}\left(\mathfrak{P}\right)$ and by Proposition~\ref{any}, we have the result.
\end{proof}
\section{The main result on generalized Lie derivations of $\mathfrak{P}$ and examples}
\begin{definition}
The set of all generalized derivations of a Lie algebra $L$ is denoted by $GenDer\left(L\right)$. It is the Lie algebra of endomorphisms $f$ in $L$ such that for every $f$ there exists another two endomorphisms $g, h$ in $L$ such that
\begin{equation}
\left[f\left(X\right),Y\right]+\left[X,h\left(Y\right)\right] = g\left[X,Y\right]\ \label{bneta}
\end{equation}
for all $X,Y\in L.$
\end{definition}
\begin{definition}
The center of a Lie algebra $L$ is the set $Z(L) = \{X\in L \mid \left[X, L\right] = \{0\}\}$.
\end{definition}
By~\cite{luk}, It is clear that $QDer\left(L\right)$ and $GenDer\left(L\right)$ are two sets of Lie algebras and satisfy $QDer\left(L\right)\subset GenDer\left(L\right)$. Moreover by assuming that a quasiderivation is a generalization of derivation and a generalized derivation is a generalization of quasiderivation, we have following result for a centerless Lie algebra $L$:
\[
ad_L\subseteq Der\left(L\right)\subseteq QDer\left(L\right)\subseteq GenDer\left(L\right)\subseteq \mathfrak{gl}\left(L\right).
\]
As we know that the center of $\mathfrak{P}$ is null by Proposition~\ref{lova}, it is interesting to find whether above inclusions verifies for $\mathfrak{P}$ or not?.

For our convenience, we denote $f\in GenDer\left(L\right)$ by $(f,h,g)$, $f\in QDer\left(L\right)$ by $(f,f,g)$ and $f\in QC\left(L\right)$ by $(f,-f,0)$. Study in~\cite{luk}, holds only for finite dimensional Lie algebras. But the following equation
\[
(f,h,g) = \left(\frac{f+h}{2},\frac{f+h}{2},g\right)+\left(\frac{f-h}{2},\frac{h-f}{2},0\right)
\]
is always holds for our $L = \mathfrak{P}.$ It further yields the following results.
\begin{proposition}
$GenDer\left(\mathfrak{P}\right) = QDer\left(\mathfrak{P}\right)+QC\left(\mathfrak{P}\right)$.\label{gene}
\end{proposition}
\begin{theorem}\label{noe}If $\mathfrak{P}$ verifies all hypotheses in the Proposition~\ref{ita}, then all generalized derivations of $\mathfrak{P}$ are of the form
\[
h+\left(L_{X+Y},L_{X+Y},L_{X+Y}+k\right)+\left(f'',f'',0\right)
\]
with $h\in QC\left(\mathfrak{P}\right)\subset\left\langle Id, Id_1, \dots, Id_n \right\rangle_\mathbb{R}$ if
\[
\left\langle x^1 \frac{\partial}{\partial x^1},\dots,x^n \frac{\partial}{\partial x^n}\right\rangle_\mathbb{R}\subset\mathfrak{P}_0
\]
and $h\in \left\langle Id\right\rangle_\mathbb{R}$ if
\[
\mathfrak{H}_0 = \mathfrak{P}_0.
\]
Where $X,Y\in\mathfrak{P}$, $k\in End\left(\mathfrak{P}\right)$ is null except on $\mathfrak{P}_0\ominus\left[\mathfrak{P},\mathfrak{P}\right]$ where $k$ is arbitrary, $f''$ is a homogeneous map of degree $1$, given in Proposition~\ref{unie}. $f''$ is null if $\mathfrak{P}_0$ contains a non diagonal linear vector field.
\end{theorem}
\begin{proof}
By corollary~\ref{bon}, Propositions~\ref{gene},~\ref{gena},~\ref{any},~\ref{uni},~\ref{unie},~\ref{ita} and~\ref{en}, we have the final result taking into account that $X\in \mathfrak{H}_0$ must be in the normalizer $\mathfrak{N}$ of $\mathfrak{P}$, because $L_X$ is stable in $\mathfrak{P}$. By Theorem~\ref{julie}, $\mathfrak{N} = \mathfrak{P}$ because all linear diagonal vector fields are in $\mathfrak{P}$ which is also separated.
\end{proof}
\begin{example}
If we take the example of Example~\ref{fs}, we have
\[
GenDer\left(\mathfrak{P}\right) = \left\langle Id, Id_1, Id_2\right\rangle_\mathbb{R}+ad_{\mathfrak{P}}+ G+ K.
\]
Where $G$ is the space of $f''$ in Proposition~\ref{unie} and $K$ is the space of quasiderivations of the form $(0, 0, g)$, where $g$ is arbitrary and $g$ vanishes except on $\left\langle x\frac{\partial}{\partial x},y\frac{\partial}{\partial y}\right\rangle_\mathbb{R}$.
\end{example}
\begin{example}
In $\mathbb{R}^2$, we consider the habitual coordinate system $(x,y)$. The Lie algebra $\mathfrak{P}$ is the infinite dimensional vector space
\[
\left\langle\frac{\partial}{\partial x},\frac{\partial}{\partial y},x\frac{\partial}{\partial x},y\frac{\partial}{\partial y},y\frac{\partial}{\partial x},y^2 \frac{\partial}{\partial x},y^3 \frac{\partial}{\partial x},\dots \right\rangle_\mathbb{R},
\]
then
\[
GenDer\left(\mathfrak{P}\right) = \left\langle Id\right\rangle_\mathbb{R}+ad_{\mathfrak{P}}+K.
\]
Where $K$ is same as in the previous example.
\end{example}

\EditInfo{January 04, 2022}{June 14, 2022}{Friedrich Wagemann}

\end{paper}